# A Stochastic Conjugate Gradient Method for Approximation of Functions

Hong Jiang[*] and Paul Wilford
Bell Laboratories, Alcatel-Lucent
700 Mountain Ave, P.O. Box 636, Murray Hill, NJ 07974-0636

**Abstract** -- A stochastic conjugate gradient method for approximation of a function is proposed. The proposed method avoids computing and storing the covariance matrix in the normal equations for the least squares solution. In addition, the method performs the conjugate gradient steps by using an inner product that is based stochastic sampling. Theoretical analysis shows that the method is convergent in probability. The method has applications in such fields as predistortion for the linearization of power amplifiers.

## 1. Introduction

This paper is concerned with the best approximation of a function by a linear combination of a set of linearly independent basis functions, e.g, by a polynomial of a certain degree.

An example of such an approximation problem arises from digital predistortion for linearization of high power amplifiers (HPA). A key step in digital predistortion is to find an approximation of the inverse function of the HPA by, e.g., a polynomial, see [1,2] and references given therein. The inverse function of HPA is not known explicitly (neither is HPA itself), but it can be observed by monitoring samples of the input and output signals to and from the HPA. The samples are usually from a waveform such as an OFDM signal, which can be modeled as a random process of a Rayleigh distribution. Furthermore, since the OFDM signal is continuously transmitted, one can always capture samples when desired. Other examples of such function approximation can be found in [3].

Therefore, the approximation problem we are interested can be characterized with the following properties. 1) The function to be approximated is not known explicitly. 2) The input and output of the function can be observed with samples. The input samples are from a random process with a certain probability density function, and its distribution can be observed, but cannot be controlled or altered. 3) There is an unlimited supply of input and output samples for observations, but there may be a limit on how many samples one can observe at one time.

Since the function to be approximated is known only through observations, the approximation of the function is computed for each set of samples from an observation. This naturally defines an iterative process in which a series of approximations are computed for a series of sample sets, with the expectation that the approximations get progressively more accurate as more sample sets are taken.

The best approximation in a vector space is normally carried out by the least squares method in which a linear combination of the basis functions is sought so that it best matches the observed output samples when evaluated at the observed input samples. The coefficients of the least squares solution satisfy the normal equations. The normal equations can be solved by an iterative method such as the conjugate gradient (CG) method. An iterative method for solving the normal equations has many advantages over a direct method such as the Cholesky decomposition. When the CG method is used to solve the normal equations, the solution to the approximation problem of the interest becomes an inner-outer loop. In the outer loop, each iteration consists of taking a set of samples from the observation and forming the normal equations with the samples. The inner loop is the CG iteration.

---

[*] Corresponding author.
  Email addresses: hong.jiang@alcatel-lucent.com (H. Jiang), paw@alcatel-lucent.com (P. Wilford)

A conjugate gradient method involving random data, such as noise, stochastic sampling, digital signal in wireless communications etc, is conventionally called a Stochastic Conjugate Gradient (SCG) method. SCG methods have been used in the literature to solve systems of linear equations from many different applications. A CG method with data from noisy measurements is considered in [4] where improvements are made to stabilize the conjugate gradients in the noisy Hessian calculation. In [5], a frame algorithm using adaptive iterative directions is used to alleviate the issue of noisy evaluation of the gradient direction. Equalization in wireless communications is considered in [6,7,8]. A nonlinear system is solved iteratively in an inner-outer loop in which a linear system of equations with the Hessian as the matrix is formed in the outer loop, and the CG iteration is used in the inner loop to solve the linear system.

In a CG or SCG method, a matrix by vector multiplication is needed. In many applications, it is possible to perform the matrix-vector multiplication directly without the explicit computation and storage of the matrix. In [9], fast curvature matrix-vector products are formed to avoid storage of the Hessian. The work in [10] has investigated approximations that can be used to efficiently perform matrix-vector multiplications when computing electric fields. The computation of the Jacobian is avoided in [10] because the product of the Jacobian with a vector can be computed by using a perturbation in the search direction for approximation of derivatives. A similar approach is used in [6,7] for Hessian. There are also many other publications, e.g., [11, 12], on improving the efficiency of conjugate gradient type methods for different applications. However, since the matrix has different structure in different applications, there is no universal matrix-vector multiplication method for avoiding the matrix computation and storage. A specific application requires a specific algorithm for efficient matrix-vector multiplication without explicit usage of the matrix.

In this work, we present a stochastic conjugate gradient for the application of approximating a function as stated at the beginning of this section. In the SCG method of this work, only one CG iteration is performed for each set of samples taken from observation. After each CG iteration, the update is used in the outer loop, and a new set of samples are taken to compute the search direction for the next update.

There are two objectives for this paper. First, a transformation is used in this work so that the conjugate gradient computations are performed in a function space instead of the traditional Euclidean vector space. Therefore, the iterations are carried out directly on the functions themselves, rather than on the coefficients of the basis functions. Because of this transformation, explicit computation and storage of the covariance matrix are avoided for this application where the matrix is traditionally computed and stored [1,2,8]. This significantly reduces the complexity in the computations.

Secondly, in this work, the inner products used in the conjugate gradient computations are approximated by a stochastic sampling method. The use of the sample evaluated inner products provides further efficiency for the calculations, and this is a distinctive feature of the SCG algorithm proposed in this paper. It is mainly because of this feature that the algorithm of this paper is called a stochastic conjugate gradient method. While it is natural to approximate the inner products by using sample averages, it is not obvious whether the SCG using only one iteration per sample set would converge. Theoretical proof will be given to show that the SCG method of this work is convergent in probability. Simulations are performed to confirm the theoretical results, and the convergence is demonstrated even for very small sample size.

This paper is organized as follows. The new stochastic conjugate gradient algorithm is introduced in Section 2. The theoretical analysis of convergence will be given in Section 3. Some implementation issues will be discussed in Section 4. In Section 5, an algorithm for multivariate functions will be presented. Simulations results will be given in Section 6. Finally, the Appendix containing proofs for the theoretical statements of Section 3 will be given at the end of the paper.

## 2. The stochastic conjugate gradient method

Let $g(x)$ be a complex valued function defined on $[0,1]$. Let $Y(t)$ be a random process with a probability density function $\rho(x) > 0, x \in [0,1]$ and $\int_0^1 \rho(x)dx = 1$. Let $Z(t) = g(Y(t))$. We assume $\rho(x)$ is known. However, the function $g(x)$, which is to be approximated, is not known explicitly, but the random process $Z(t)$ can be observed, and its samples can be taken as desired.

Let $\{\phi_0(x),...,\phi_{M-1}(x)\}$ be a set of linearly independent complex valued functions of real variables defined on the interval $[0,1]$. Let $\phi = [\phi_0(x),...,\phi_{M-1}(x)]^T$. For example, $\{\phi_0(x), \phi_1(x)...,\phi_{M-1}(x)\} = \{1, x..., x^{M-1}\}$ forms a basis for the polynomials of degree less than $M$. Let $\mathbf{E}(y)$ denote the expected value of the random variable $y$, and $y^*$ the complex conjugate of $y$. The problem that we are interested can be stated as follows.

**Problem 1**

Find $\bar{u}(x) = \sum_{i=0}^{M-1} \bar{u}_i \phi_i(x)$, such that

$$\mathbf{E}\left((g(Y) - \bar{u}(Y))^*(g(Y) - \bar{u}(Y))\right) = \mathbf{E}\left(\left(Z - \sum_{i=0}^{M-1} \bar{u}_i \phi_i(Y)\right)^* \left(Z - \sum_{i=0}^{M-1} \bar{u}_i \phi_i(Y)\right)\right)$$
$$= \min\{\mathbf{E}\left((Z - u(Y))^*(Z - u(Y))\right) | u \in \mathbf{V}\}, \quad (2.1)$$
$$\mathbf{V} = \text{span}\{\phi_0(x), \phi_1(x)..., \phi_{M-1}(x)\}.$$

The solution to Problem 1 can be readily obtained through a minimization process. Define the inner product $\langle \cdot, \cdot \rangle$ of complex valued functions on the interval $[0,1]$ by

$$\langle u, v \rangle = \mathbf{E}(u^* v) = \int_0^1 \rho(x) u^*(x) v(x) dx. \quad (2.2)$$

Define a functional of functions define on $[0,1]$ by
$$J(u) = \langle u - g, u - g \rangle. \quad (2.3)$$

Then the solution $\bar{u}$ to Problem 1 is equivalent to the solution to the following minimization problem
$$J(\bar{u}) = \min\{J(u) | u \in \mathbf{V}\}. \quad (2.4)$$

It is well known that solving (2.4) is equivalent to solving the normal equations
$$Au = b, \quad (2.5)$$
where $A = A(\phi)$ is the $M \times M$ covariance matrix, and $u = [u_0,...,u_{M-1}]^T$ and $b = [b_0,...,b_{M-1}]^T$ are complex $M$-tuplets. Their components are given by
$$a_{ij} = \langle \phi_i, \phi_j \rangle, \quad b_i = \langle \phi_i, g \rangle. \quad (2.6)$$

Therefore, the solution to Problem 1 is completely determined if the inner products in (2.6) are known.

The solution $\bar{u} = [\bar{u}_0,...,\bar{u}_{M-1}]^T$ to the normal equations (2.5) can be found iteratively by using the conjugate gradient method. After the coefficients are computed, the approximating function $\bar{u}(x)$ can be obtained by $\bar{u}(x) = \sum_{i=0}^{M-1} \bar{u}_i \phi_i(x)$. The details of the CG method can be found, e.g., in [13]. Instead of computing the coefficients $[\bar{u}_0,...,\bar{u}_{M-1}]^T$ first, and then computing the approximating function $\bar{u}(x)$ using the coefficients, it is possible to compute the approximating function $\bar{u}(x)$ directly without computing the coefficients. For this purpose, the conjugate gradient method for the solution of (2.4) can be derived with a minimization process working on the function space $\mathbf{V} = \text{span}\{\phi_0,...,\phi_{M-1}\}$ directly. The resulting algorithm is given as follows.

**Algorithm CG**:
$u^0 = 0; v^0 = 0;$
loop for $k = 1, 2, ..., M$

$$r^{k-1} = \sum_{j=0}^{M-1} \langle \phi_j, g - u^{k-1} \rangle \phi_j \qquad (2.7)$$

terminate when $\langle r^{k-1}, r^{k-1} \rangle = 0$

$$\beta_k = \langle r^{k-1}, r^{k-1} \rangle / \langle r^{k-2}, r^{k-2} \rangle, \quad k > 1 \qquad (2.8)$$

$$v^k = r^{k-1} + \beta_k v^{k-1} \qquad (2.9)$$

$$\alpha_k = \langle r^k, g - u^{k-1} \rangle / \langle v^k, v^k \rangle \qquad (2.10)$$

$$u^k = u^{k-1} + \alpha_k v^k \qquad (2.11)$$

end loop

After the termination of Algorithm CG, $u^k$ is the solution to Problem 1. As a convention in this paper, unless otherwise stated, we use the superscript for the index of CG iteration if the variable is a vector or a function, e.g., $u^k$, so that the subscript can be used as the index for components of the vector. We will use the subscript for the index of CG iteration if the variable is a scalar, e.g., $\alpha_k$.

Algorithm CG above computes the functions $r^{k-1}$, $v^k$ and $u^k$ in the function space **V**, as opposed to the classic CG method which computes complex *M*-tuplets of the Euclidean space $C^M$. In Algorithm CG, equation (2.7) computes the residual $r^{k-1}$ from the previous approximation $u^{k-1}$. The new search direction $v^k$ is computed in (2.9). Initially, the search direction is equal to the residual. Subsequently, the new search direction is orthogonal to the previous search direction in the inner product $\langle \cdot, \cdot \rangle$ defined in (2.2), which is conjugate-orthogonal in the Euclidean inner product of the *M*-tuplets. The new approximation $u^k$ is computed in (2.11) to minimize the functional $J(u)$ in the search direction of $v^k$. The algorithm terminates in no more than $M$ iterations in the absence of round-off errors. Function $u^k$ after the termination of the iteration is the solution to Problem 1. Note that Algorithm CG only requires the computation of the residual in (2.7), it does not require the computation or storage of the covariance matrix $A$.

In practice, Algorithm CG cannot be realized if function $g(x)$ to be approximated is not known explicitly. This is because the computations in (2.7) and (2.10) involve the inner product of $g(x)$. Therefore, the inner products in (2.7) and (2.10) must be approximated in a way that they become computable.

Let $\{z_0, ..., z_{N-1}\}$, and $\{y_0, ..., y_{N-1}\}$, $n = 0, ..., N-1$, be sets of samples taken from the random processes $Z(t)$ and $Y(t)$, respectively, so that $z_n = g(y_n)$. The inner product of any two functions is given by (2.2). Evaluating the functions $u(x), v(x)$ at the samples $\{y_0, ..., y_{N-1}\}$, and defining

$$u(y) = [u(y_0), ..., u(y_{N-1})]^T, \text{ and } v(y) = [v(y_0), ..., v(y_{N-1})]^T, \qquad (2.12)$$

we have

$$\langle u, v \rangle = \mathbf{E}(u(Y)^* v(Y)) = \lim_{N \to \infty} \frac{1}{N} u(y)^H v(y).$$

Therefore, the inner product $\langle u, v \rangle$ may be approximated by

$$\langle u, v \rangle \approx \tfrac{1}{N} u(y)^H v(y) = \tfrac{1}{N} \sum_{n=0}^{N-1} u(y_n)^* v(y_n).$$

For any two functions $u, v$ defined on $[0,1]$, we define

$$\langle u, v \rangle_s \overset{\Delta}{=} \tfrac{1}{N} \sum_{n=0}^{N-1} u(y_n)^* v(y_n). \tag{2.13}$$

Strictly speaking, $\langle u, v \rangle_s$ is not an inner product because the support of $u(x), v(x)$ may not intersect the sample set $\{y_0, ..., y_{N-1}\}$. In addition, the value of $\langle u, v \rangle_s$ depends on a particular instance of samples. However, it is reasonable to assume that when $u(x)$ is continuous and when the sample size $N$ is large enough, at least one sample will be in the support of $u(x)$, and hence $\langle \cdot, \cdot \rangle_s$ of (2.13) defines an inner product. This inner product differs from sample set to sample set, i.e., in general, $\langle u, v \rangle_{s_1} \neq \langle u, v \rangle_{s_2}$, where $s_1, s_2$ represent two different sample sets.

Some modifications need to be made to Algorithm CG for it to be practical. First, the inner product $\langle \cdot, \cdot \rangle$ will be approximated by $\langle \cdot, \cdot \rangle_s$ of (2.13). Secondly, there is a need to restart the iteration process, i.e., reset the search direction $v^k$ to residual $r^{k-1}$, from time to time. The reason is that the search directions may no longer be orthogonal to each other because the inner products $\langle \cdot, \cdot \rangle_s$ are different from iteration to iteration. Restarting a Krylov subspace method is a common practice, e.g., with GMRES [14]. A simple strategy is to restart the process after a predetermined number of iterations have been performed. Thirdly, the calculation of $\beta_k$ in (2.8) needs to be revised to guarantee orthogonality between $v^k$ and $v^{k-1}$ when a different inner product is used at a different iteration. The following algorithm is derived from Algorithm CG with these modifications, and the notation $\langle \cdot, \cdot \rangle_k$ is used to signify that the inner product is computed as (2.13) with the sample set $k$.

**Algorithm SCG**
At start:
    Given threshold $\varepsilon > 0$
    Determine a strategy to reset the search direction at least once every $M$ iterations
    $u^0 = 0; v^0 = 0;$
loop for $k = 1, 2, ...$
    take sample sets $\{y_0, ..., y_{N-1}\}, \{z_0, ..., z_{N-1}\}$ such that $z_n = g(y_n), n = 0, ..., N-1$

$$r^{k-1} = \sum_{j=0}^{M-1} \langle \phi_j, g - u^{k-1} \rangle_k \phi_j \tag{2.14}$$

    if at start, or at reset, then

$$\beta_k = 0 \tag{2.15}$$

    else

$$\beta_k = -\langle r^{k-1}, v^{k-1} \rangle_k / \langle v^{k-1}, v^{k-1} \rangle_k \tag{2.16}$$

    end if
    $v^k = r^{k-1} + \beta_k v^{k-1} \tag{2.17}$

    if $\langle v^k, v^k \rangle_k < \varepsilon$, then
        reset at next iteration and go to next iteration
    end if

$$\alpha_k = \langle r^k, g - u^{k-1} \rangle_k / \langle v^k, v^k \rangle_k \tag{2.18}$$

$$u^k = u^{k-1} + \alpha_k v^k \tag{2.19}$$

end loop

At each iteration $k$, $u^k(x)$ is an approximation of the solution to Problem 1. Algorithm SCG has no stopping criteria; it provides approximation of $g(x)$ continuously. At each iteration, a number of samples are taken. The sample size, $N$, may be different from iteration to iteration. It is also possible that the same samples are used from the previous iteration. This corresponds to performing more than one iteration of SCG with the same inner product. If the same sets of samples are kept unchanged for $M$ iterations, the residual is guaranteed to be zero in no more than $M$ iterations in the absence of round-off errors, although the resulting function $u^k$ may still not be the solution to Problem 1, because the resulting function $u^k$ minimizes $J_k(u) = \langle u - g, u - g \rangle_k$, but not necessarily $J(u) = \langle u - g, u - g \rangle$.

The formulas given in Algorithm SCG are for the convenience of presentation, and they can be rewritten for more computational efficiency. For example, the coefficient of $\phi_j$ in (2.14) can be computed as

$$\gamma_j = \tfrac{1}{N} \sum_{n=0}^{N-1} \phi_j(y_n)^* (z_n - u^{k-1}(y_n)) \quad j = 0,\ldots,M-1. \tag{2.20}$$

And the numerator of (2.18) can be simplified as

$$\langle r^{k-1}, g - u^{k-1} \rangle_k = \sum_{j=0}^{M-1} \gamma_j^* \gamma_j. \tag{2.21}$$

Also,

$$\langle v^{k-1}, v^{k-1} \rangle_k = \tfrac{1}{N} \sum_{n=0}^{N-1} v^{k-1}(y_n)^* v^{k-1}(y_n), \text{ and } \langle v^k, v^k \rangle_k = \tfrac{1}{N} \sum_{n=0}^{N-1} v^k(y_n)^* v^k(y_n). \tag{2.22}$$

If the reset (2.15) is performed at every iteration, Algorithm SCG becomes a stochastic steepest decent method because the search direction $v^k$ is always equal to the residual $r^{k-1}$ which is the direction of the steepest decent.

## 3. The convergence analysis

In this section, we discuss convergence properties of Algorithm SCG. While it is natural to approximate the inner product (2.2) by the sample average (2.13), the question also arises as to whether Algorithm SCG actually converges. We will show theoretically that the approximation computed from Algorithm SCG converges to the solution of Problem 1 in probability.

We start with some necessary definitions. Let $\bar{u} \in \mathbf{V} = \text{span}\{\phi_0,\ldots,\phi_{M-1}\}$ be the solution to Problem 1. We define a functional on $\mathbf{V}$ as

$$H(u) = \langle u - \bar{u}, u - \bar{u} \rangle. \tag{3.1}$$

The error in the approximation after $k$ iterations of Algorithm SCG is $u^k - \bar{u}$, and $H(u^k)$ is therefore the square of the norm of the error, or mean square error. Therefore, Algorithm SCG is convergent if and only if $\lim_{k \to \infty} H(u^k) = 0$. We denote the condition number[1] of the covariance matrix $A$ of (2.6) by $c = \text{cond}(A)$.

---

[1] The condition number is traditionally denoted by the Greek letter $\kappa$, but it is too easily confused with the iteration index $k$.

Next, we define the counterparts of these in the sample evaluated inner product $\langle \cdot, \cdot \rangle_k$ which is defined in (2.13) with the samples taken at iteration $k$ of Algorithm SCG. First, we define the counterparts of $J(u)$ and $H(u)$ as

$$J_k(u) = \langle g - u, g - u \rangle_k, H_k(u) = \langle u - \bar{u}^k, u - \bar{u}^k \rangle_k. \tag{3.2}$$

In (3.2), $\bar{u}^k$ is the solution to the minimization problem

$$J_k(\bar{u}^k) = \min\{J_k(u) \mid u \in \mathbf{V}\}. \tag{3.3}$$

It can be shown that for every $u \in \mathbf{V}$

$$H_k(u) = J_k(u) - \langle g - \bar{u}^k, g - \bar{u}^k \rangle_k. \tag{3.4}$$

For a given $k$, the last term in (3.4) is a constant, independent of $u$. At each SCG iteration $k$, a covariance matrix $A^k$ can be defined similarly to $A$ of (2.6) except that $\langle \cdot, \cdot \rangle$ is replaced by $\langle \cdot, \cdot \rangle_k$. Matrix $A^k$ is not needed in Algorithm SCG, but it is convenient in the discussion of convergence properties. The condition number of $A^k$ is denote by $c_k = \text{cond}(A^k)$. We will assume that the sample size taken at each iteration of Algorithm SCG is large enough so that $\langle \cdot, \cdot \rangle_k$ approximates $\langle \cdot, \cdot \rangle$, and the condition numbers of $A^k$ is in the same order as that of $A$. More precisely, we make the following assumptions.

**Assumptions**

We assume that the following properties hold for $\langle \cdot, \cdot \rangle_k$.

1) $\langle \cdot, \cdot \rangle_k$ is an inner product, i.e., for any function $u$,

$$\langle u, u \rangle_k = 0 \text{ implies } u = 0. \tag{3.5}$$

2) For $k \geq 1$ and $u \neq 0$, the random variables

$$\omega = \omega(k, u) = \frac{\langle u, u \rangle_k}{\langle u, u \rangle}, \tag{3.6}$$

are independent and identically distributed. Their logarithms have a finite mean and variance, and they are given by

$$\mathbf{E}(\ln(\omega(k,u))) = \mu,$$
$$\mathbf{Var}(\ln(\omega(k,u))) = \frac{\sigma^2}{2}. \tag{3.7}$$

The values $\mu, \sigma^2$ are related to the number of samples. If the sample size is large enough so that $\langle \cdot, \cdot \rangle_k$ well approximates $\langle \cdot, \cdot \rangle$, then $\mu, \sigma^2$ have small values because $\omega$ of (3.6) is close to 1.

3) The condition numbers $c_k$ of $A^k$ have an upper bound. The solutions $\bar{u}^k$ of (3.3) have an upper bound. That is, there exist $c_0: c_0 \geq c \geq 1$, and $d_0 > 0$ such that

$$c_k \leq c_0, \langle \bar{u}^k, \bar{u}^k \rangle < d_0^2, \text{ for all } k \geq 1. \tag{3.8}$$

We are now at a position to state some properties of Algorithm SCG. All proofs are postponed to Appendix.

**Lemma 1.** Let $u^k$ be the approximation after $k \geq 1$ iterations of Algorithm SCG. Then, there exist a sequence $\delta_k$ with

$$0 \leq \delta_k \leq 1, \tag{3.9}$$

such that
$$H_k(u^k) = \delta_k H_k(u^{k-1}), \quad k = 1, 2, \ldots. \tag{3.10}$$
Furthermore, there exists a $p$ with $0 < p < 1$ such that
$$\prod_{j=1}^{k} \delta_j \leq \left(1 - \frac{1}{c_0}\right)^{pk}, \text{ for large } k, \tag{3.11}$$
where $c_0$ is the bound on the condition numbers given in (3.8).
**Proof.** See Appendix.

Equations (3.9) and (3.10) show that $H_k(u^k)$ is reduced at each iteration, while (3.11) provides an estimate on the rate of reduction. It is worthwhile to point out that the upper bound on the right hand side (RHS) of (3.11) is too pessimistic. The left hand side (LHS) of (3.11) is expected to be much smaller than RHS. The base, $1 - c_0^{-1}$, is too pessimistically large. As is shown in the proof of Lemma 1, the upper bound in (3.11) is derived essentially based on the steepest decent method, but we expect Algorithm SCG to perform better. Also, the factor in the exponent, $p$, could be very close to 1. Regardless, (3.11) is all we need to establish the convergence of Algorithm SCG, although in reality, the convergence can be much faster than predicted by it. Lemma 1 alone, however, does not imply that $u^k$ will converge to the solution of Problem 1, because the reduction is measured according to $\langle \cdot, \cdot \rangle_k$, which varies from iteration to iteration. Lemma 1 does not guarantee that $H(u^k)$ is also reduced. The next Lemma is to address the stochastic nature of $\langle \cdot, \cdot \rangle_k$.

**Lemma 2.** If $H(u^k) \neq 0$ for all $k = 1, 2, \ldots$, then there exist a sequence $\eta_k$ and a random process $\theta_k$ such that
$$H(u^k) = \eta_k \theta_k H(u^0). \tag{3.12}$$
Furthermore, there exists a $p_0$ with $0 < p_0 < 1$, so that $\eta_k$ satisfies
$$\eta_k \leq \left(1 - \frac{1}{c_0}\right)^{p_0 k}, \text{ for large } k. \tag{3.13}$$
For large $k$, the random variable $\theta_k$ has a log-normal distribution, and its probability density function is given by
$$f_{\theta_k}(t) = f_{\text{LogN}}(t; 0, k\sigma^2) = \frac{1}{t\sqrt{k}\sigma\sqrt{2\pi}} e^{-\frac{(\ln t)^2}{2k\sigma^2}}, \quad t > 0, \tag{3.14}$$
where $\sigma^2$ is given in (3.7).
**Proof.** See Appendix.

The two factors on the RHS of (3.12) characterize the asymptotic behaviors of the two processes involved in the SCG method. The first factor $\eta_k$ is the result of the conjugate gradient iterative process in which the mean square error, $H_k(u^k)$, is reduced by that amount after $k$ iterations. The second factor $\theta_k$ captures the impact of the stochastic process of replacing the inner product $\langle \cdot, \cdot \rangle$ by its sample evaluated version $\langle \cdot, \cdot \rangle_k$. Although the first factor gets progressively small as the iteration number is increased, the second factor $\theta_k$ may get unboundedly large as the iteration is increased. In fact, as $k \to +\infty$, $\mathbf{E}(\theta_k) \to +\infty$. The two factors, $\eta_k$ and $\theta_k$, counteract against each other, but the rate of decay of $\eta_k$ is faster than the rate of growth of $\theta_k$, so Algorithm SCG converges, which is stated in the following theorem.

**Theorem 1.** The computed function $u^k$ from Algorithm SCG converges to the solution of Problem 1 in probability. More precisely, we have
1) either $H(u^k) = 0$ for some $k$, in which case $u^k = \bar{u}$ is the solution to Problem 1, or
2) $\lim_{k \to +\infty} H(u^k) = 0$ in probability, i.e., for every $\varepsilon > 0$,

$$\lim_{k \to +\infty} \Pr\left(\left|H(u^k)\right| < \varepsilon\right) = 1. \tag{3.15}$$

In (3.15), $\Pr(\cdot)$ is the probability.
**Proof.** See Appendix.

As shown in the proof of Theorem 1, there are two numbers that determine the convergence rate of Algorithm SCG. The first is $1 - c_0^{-1}$; the smaller this number is, the faster the convergence is. The parameter $c_0$ is related to the condition number $c$ of the covariance matrix $A$, and a reduced $c$ speeds up convergence. Therefore, it is important to properly choose the basis functions $\{\phi_0, \ldots, \phi_{M-1}\}$ to reduce the condition number of the resulting covariance matrix $A$, even though the matrix does not explicitly appear in Algorithm SCG.

The second number that affects the convergence rate is the variance $\sigma^2$ of the logarithm of the random variables defined in (3.6) and (3.7). The smaller $\sigma^2$ is, the faster the convergence is. This number is determined by the sample sets used in the evaluation of the inner products $\langle \cdot, \cdot \rangle_k$. A large sample size implies a small variance $\sigma^2$, and hence a fast convergence rate for Algorithm SCG. Also, when $\sigma^2$ is small, the variance in the approximation $u^k$ is small, and the convergence will be smooth. In fact, if $c_0$ and $\sigma^2$ are small enough so that

$$\left(1 - \tfrac{1}{c_0}\right)^{p_0} e^{\sigma^2/2} < 1, \tag{3.16}$$

then $u^k$ converges to $\bar{u}$ in mean square, i.e.,

$$\lim_{k \to +\infty} \mathbf{E}\bigl(H(u^k)\bigr) \le H(u^0) \lim_{k \to +\infty} \left(\left(1 - \tfrac{1}{c_0}\right)^{p_0 k} \mathbf{E}(\theta_k)\right) = H(u^0) \lim_{k \to +\infty} \left(\left(1 - \tfrac{1}{c_0}\right)^{p_0} e^{\sigma^2/2}\right)^k = 0. \tag{3.17}$$

Equation (3.17) is a directly consequence of Lemma 2 because $\mathbf{E}(\theta_k) = e^{\sigma^2/2}$. Equation (3.17) represents a stronger result than Theorem 1 because of the additional assumption (3.16).

## 4. Implementation considerations

### 4.1 Orthogonal basis functions

Although Theorem 1 guarantees convergence of Algorithm SCG, the convergence rate depends on condition numbers of the covariance matrices. This is clear from equations (3.12) and (3.13). The rate of convergence is bounded by $\eta_k$ as given in (3.12), and in turn, the bound of $\eta_k$ depends on $c_0$, as is shown by (3.13), where $c_0$ is the bound of condition numbers $c_k = \text{cond}(A^k)$ given in (3.8). The larger the condition numbers are, the slower the convergence rate is. In many applications, the covariance matrices may be ill-conditioned, resulting in slow convergence. This problem can be alleviated by using an orthogonal basis.

Orthogonal basis functions help reducing the condition number of the covariance matrix, and hence improving convergence of Algorithm SCG. The basis functions $\{\phi_0,...,\phi_{M-1}\}$ may be orthogonalized by, for example, the Gram-Schmidt process, in which an orthonormal basis $\{\psi_0,...,\psi_{M-1}\}$ is obtained with the property

$$\psi_k(x) \in \text{span}\{\phi_0(x),...,\phi_k(x)\}, k = 0,...,M-1,$$

$$\langle \psi_i(x), \psi_j(x) \rangle = \begin{cases} 0 & i \neq j \\ 1 & i = j \end{cases}.$$

For polynomials with the basis functions $\{\phi_0(x), \phi_1(x)..., \phi_{M-1}(x)\} = \{1, x,..., x^{M-1}\}$, an orthonormal basis $\{\psi_0,...,\psi_{M-1}\}$ can be constructed simply by using a three term recursion. Orthogonal polynomials are also considered in [15] for the application of predistortion.

With an orthonormal basis, the covariance matrix has condition number $\text{cond}(A(\psi)) = 1$. However, this does not imply that the solution of the least squares problem can be found after one iteration of SCG. The reason is that the basis functions that are orthogonal in $\langle \cdot, \cdot \rangle$ may no longer be orthogonal in $\langle \cdot, \cdot \rangle_k$, and therefore, the condition number of $A^k$ is in general greater than one. Nevertheless, we expect the condition number of $A^k$ to be smaller when the basis functions are orthogonal.

### 4.2 Estimate for the probability density function

The probability density function $\rho(x)$ of the random process $Y(t)$ is needed in the orthogonalization process. In general, a histogram or the kernel density estimation may be used to estimate the probability density function. In the application of predistortion, the signal is usually an OFDM signal. The amplitude of an OFDM signal has a distribution close to a Rayleigh distribution, and therefore, we may choose the probability density function as $\rho(x) = xe^{-x^2/(2\sigma^2)}/\sigma^2$. The parameter $\sigma$ may be estimated from the observed samples by $\hat{\sigma} = \sqrt{\frac{1}{2N}\sum_{n=0}^{N-1} y_n^2}$.

### 4.3 Lookup table implementation

In many applications, the computed function $u^k(x)$ from Algorithm SCG is used in further processing in an outer loop. In hardware implementation such as on an FPGA, the evaluation of a function is best accomplished by the use of a look-up table (LUT). A LUT of a function is a vector whose index corresponds to a quantized value of the independent variable of the function. Therefore, the probability density function, function $u^k(x)$, the basis functions $\{\phi_0,...,\phi_{M-1}\}$, and other functions in the SCG algorithm can all be represented by LUTs. Quantize $[0,1]$ into $B$ levels, and let $x_j = j/B, j = 0,...,B-1$. Then a look-up table representation can be defined as $LUTu(j) = u(x_j)$, for $j = 0,...,B-1$.

### 4.4 The complexity of SCG

The cost of Algorithm SCG is determined by its complexity. The complexity of Algorithm SCG can be analyzed by assuming that the functions involved are implemented by LUTs of size $B$. The majority of operations are from the computation of residual in (2.14) in which the LUT values are evaluated at $N$ samples, and then $M$ inner products $\langle \cdot, \cdot \rangle_k$ are computed, for a total of $MN$ operations. Direction $v^k$ is formed with $M$ basis functions to give $MB$ operations. The rest of the operations require two inner products $\langle \cdot, \cdot \rangle_k$, and 2 function updates. Therefore, the complexity for Algorithm SCG per iteration is $O(MN + MB + N + B)$.

This can be compared with the traditional CG algorithm, which is similar to Algorithm CG with the exceptions that 1) the residual in (2.7) is computed with a matrix by vector multiplication, and 2) each inner product is computed by vector dot product of $M$-tuplets. When the matrix is computed at each iteration, the computation requires $M(M+1)/2$ inner products, each of which requires evaluation at $N$ samples, for a total of $O(M^2 N)$ operations. The residual computation needs $MN$ operations. The rest is similar to Algorithm SCG. Therefore, the complexity per iteration for the traditional CG algorithm is $O(M^2 N + MN + MB + B)$, which is higher than the complexity of SCG.

## 5. Multivariate functions

Although theoretical treatment of multivariate functions is very close to that of functions of one variable, there are significant practical issues in multivariate functions that warrant more discussions. In this section, the superscript will be used as the index of dimensions, and will no longer be used as the iteration index. The omission of the iteration index will not cause confusions because the algorithm specifies the computations within one iteration to compute the updates for the next iteration.

For multivariate functions, even if the function to be approximated, $g(x^1, ..., x^Q)$, is known, it is no longer feasible to compute the inner product $\langle \cdot, \cdot \rangle$. Its computation not only requires multivariate integrals, but also requires knowledge of the joint probability density function $\rho(Y^1, ..., Y^Q)$, both of which may not be practical due to the curse of dimensionality. The inner product must be approximated by sampling as defined by

$$\langle u, v \rangle_s \stackrel{\Delta}{=} \tfrac{1}{N} \sum_{n=0}^{N-1} u(y_n^1, ..., y_n^Q)^* v(y_n^1, ..., y_n^Q).$$

An additive separable function is a special form of multivariate function which is a sum of functions of one variable. Additive separable functions arise in many applications. In [3], the approximation of a multivariate function is reduced to a series of problems of finding the best additive separable function approximation. In the application of predistortion for HPA with memory, it is assumed that the inverse of the HPA can be approximated by a memory polynomial of the form

$$z_n = P(x_n, ..., x_{n-Q+1}) = \sum_{q=1}^{Q} x_{n-q+1} P_q(|x_{n-q+1}|). \tag{5.1}$$

In (5.1), $x_n$ are complex samples from the transmitted signal. $P_q(\cdot)$ is a polynomial of a degree less than $M$. (5.1) is said to be a memory polynomial because $z_n$ depends on not only $x_n$ but also its past states $x_{n-1}, ..., x_{n-Q+1}$. The delays in samples are necessary to account for memory effects of a HPA with memory. The objective is, therefore, to find polynomials $P_q(\cdot)$, $q = 1, ..., Q$, so that the pair of sample sets $\{x_n\}$, $\{z_n\}$ from (5.1) best match the sample sets from observing the output and input signals of the HPA. In this context, the inverse of HPA is regarded to be a multivariate function in which a dimension is the current state or a state in the past, see [1,2] for more details.

Let $g(x^1,...,x^Q)$ be a complex valued function of $Q$ complex variables. Let $Y^1,...,Y^Q$ be complex valued random processes. Each of $|Y^q|$ has the same probability density function $\rho(x)$. Let $\mathbf{Y} = [Y^1,...,Y^Q]$. Define the random variable $Z = g(Y^1,...,Y^Q)$. Let $\{\psi_0(|x|),...,\psi_{M-1}(|x|)\}$ be a set of linearly independent complex valued functions of one real variable defined on the interval $[0,1]$. Let $\{\tau_1(x),...,\tau_Q(x)\}$ be a given set of complex valued functions of one complex variable. A generalization of the memory polynomial predistorter (5.1) is given by

$$z_n = P(x_n,...,x_{n-Q+1}) = \sum_{q=1}^{Q} \tau_q(x_{n-q+1}) \sum_{i=0}^{M-1} u_i^q \psi_i(|x_{n-q+1}|). \tag{5.2}$$

The memory polynomial predistorter in (5.1) becomes a special case of (5.2) with $\tau_q(x) = x, q = 1,...,Q$ and $\psi_i(x) = x^i$, $i = 0,...,M-1$. With these definitions, finding the best predistorter (5.2) is tantamount to solving the best approximation problem stated as follows.

**Problem 2.**

Find, $\bar{u}^q(|x|) = \sum_{i=0}^{M-1} \bar{u}_i^q \psi_i(|x|), q = 1,...,Q$ such that

$$\mathbf{E}\left(\left(g(\mathbf{Y}) - \sum_{q=1}^{Q} \tau_q(Y^q)\bar{u}^q(|Y^q|)\right)^* \cdot \left(g(\mathbf{Y}) - \sum_{q=1}^{Q} \tau_q(Y^q)\bar{u}^q(|Y^q|)\right)\right)$$

$$= \min_{v^q \in \mathbf{V}} \mathbf{E}\left(\left(Z - \sum_{q=1}^{Q} \tau_q(Y^q)v^q(|Y^q|)\right)^* \cdot \left(Z - \sum_{q=1}^{Q} \tau_q(Y^q)v^q(|Y^q|)\right)\right),$$

$$\mathbf{V} = \mathrm{span}\{\psi_0,...,\psi_{M-1}\}.$$

The solutions $\bar{u}^q(|x|), q = 1,...,Q$ to Problem 2 can be computed using the following algorithm, in which the superscript is now the index of the dimensions in a multivariate function, rather than the index of iterations as in Algorithm SCG. Since no iteration index is used, the same variable will be used both before and after an iteration, but the notation "$\leftarrow$" will be used to indicate that the variable has been assigned to a new value.

**Algorithm SCG_MUL**

At start:

    Given $\varepsilon > 0$, and basis functions $\{\psi_0(|x|),...,\psi_{M-1}(|x|)\}$

    Determine a strategy to reset the search direction at least once every $QM$ iterations.

    $u^q(x) = 0$, for $q = 1,...,Q$

loop:

    Take sample sets $\{y_0^q,...,y_{N-1}^q\}, q = 1,...,Q$, and $\{z_0,...,z_{N-1}\}$ such that $z_n = g(y_n^1,...,y_n^Q), n = 0,...,N-1$

$$\gamma_i^q = \tfrac{1}{N} \sum_{n=0}^{N-1} \left((\tau_q(y_n^q)\psi_i(|y_n^q|))^* \cdot (z_n - \sum_{q=1}^{Q} \tau_q(y_n^q)u^q(|y_n^q|))\right), q = 1,...,Q, i = 0,...,M-1$$

$$r^q = [\gamma_0^q,...,\gamma_{M-1}^q]^T, q = 1,...,Q$$

    if at start or reset

        $v^q(|x|) \leftarrow \sum_{i=0}^{M-1} \gamma_i^q \psi_i(|x|), q = 1,...,Q$

    else

$$\beta = \left(\sum_{q=1}^{Q} (r^q)^H r^q\right)/\omega,$$
$$\omega \leftarrow \sum_{q=1}^{Q} (r^q)^H r^q,$$
$$v^q(|x|) \leftarrow \sum_{i=0}^{M-1} \gamma_i^q \psi_i(|x|) + \beta v^q(|x|), q = 1,...,Q$$
  end if
  if $\omega < \varepsilon$
    reset at next iteration and go to next iteration
  else
$$\lambda = \sum_{n=0}^{N-1} \left|\sum_{q=1}^{Q} \tau_q(y_n^q) v^q(|y_n^q|)\right|^2,$$
$$\alpha = \frac{\omega}{\lambda},$$
$$u^q(|x|) \leftarrow u^q(|x|) + \alpha v^q(|x|), q = 1,...,Q$$
  end
end loop

Although the same number of basis functions, and the same basis functions are used in Algorithm SCG_MUL for each dimension $q = 1,...,Q$, it is done purely for convenience. The number of basis functions, and the basis functions themselves can be made different for different dimension $q$.

It is also advantageous to use orthonormal basis functions to reduce the condition number of the covariance matrix. For Problem 2, orthonormal basis functions are defined as

$$\langle \tau_q(x)\psi_i(|x|), \tau_q(x)\psi_j(|x|)\rangle = \begin{cases} 0 & i \neq j, q = 1,...,Q \\ 1 & i = j, q = 1,...,Q \end{cases}$$

Even if $\{\psi_0,...,\psi_{M-1}\}$ is an orthonormal basis, the condition number of the covariance matrix for the multivariate Problem 2 is still larger than one in general because the covariance matrix may not be diagonal as the dimensional variables are not necessarily independent random variables. However, experiments show that using the orthonormal basis functions for each dimension will reduce the conditional number of the covariance matrix for the multivariate problem.

## 6. Simulations

### 6.1 Function of one variable

In the first example, Algorithm SCG is applied to computing an approximation of a function of one variable given by
$$g(y) = \sin(2\pi y), y \in [0,1].$$
The basis functions $\{\psi_0,...,\psi_{M-1}\}$ are orthogonal polynomials of degree less than $M$ with respect to the uniform weight function on $[0,1]$. At each iteration, $N$ samples are taken from a Rayleigh distribution, i.e., the samples are computed by $y_n = \sqrt{(x_n^1)^2 + (x_n^2)^2}$, where $x_n^1, x_n^2$ are random numbers from the normal distribution with mean=0.25 and variance=0.0625. The values of $y_n$ that are outside of $[0,1]$ are ignored. Therefore, the uniform weight function is not the same as the probability density function of the samples. The basis functions are intentionally chosen to be not orthogonal in the inner product $\langle \cdot, \cdot \rangle$ defined by the Rayleigh distribution, but at the same time, the resulting covariance matrix still has a modest condition number. A LUT of $B$ entries is used to represent each function involved in Algorithm SCG. The algorithm is reset after every $M$ iterations. Only one iteration is performed for each set of samples captured. At each iteration, the mean square error, $H(u^k)$, is obtained by computing $u^k(y) - \sin(2\pi y)$ at 1000 evenly spaced points in the interval $[0,1]$, and therefore, the measurement $H(u^k)$ is independent of the sample set. We report the values of $H(u^k)$, which is computed with the uniform weight function in the inner product, as a function of the iteration number $k$. Four simulations are performed in which a different sample size $N$ is used for each simulation run. The parameters used in the simulations are given as $M = 10, B = 2^{16}, N = 1, 50, 100, 500$. In other words, the sine function is approximated by a polynomial of degree 9. The simulation results for first 100 iterations are shown in Figure 1.

For the simulation with sample size $N = 500$, the mean square error $H(u^k)$ reaches to the level of about $4.0e - 10$ in less than 30 iterations, and after that it remains at that level almost as a straight horizontal line. The error cannot be reduced further for two reasons: 1) the computed function $u^k(y)$ is represented by a LUT of 16 bits, which limits the resolution of LUT and hence the accuracy of $u^k(y)$, and 2) the mean square error $H(u^k)$ is bounded below by the square of the distance between the function $\sin(2\pi y)$ and the linear space of polynomials of degrees less than $M = 10$. To demonstrate this further, the curve for $N = 1000$ (not shown) after convergence is almost indistinguishable from that of $N = 500$.

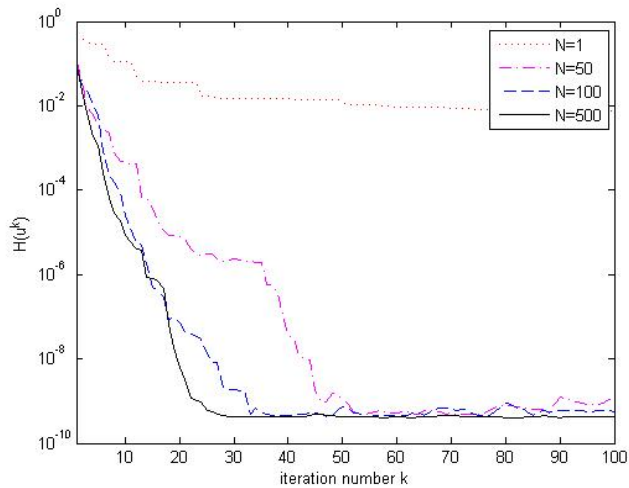

Figure 1. Convergence of mean square error: from top to bottom, sample size N=1, 50,100, 500

The errors for $N = 100$ and $N = 50$ also reach to the level of $4.0e - 10$, in no more than 50 iterations, but for $N = 1$, it takes much longer for the error to reach the same level, as shown in Figure 2. Nevertheless, it is quite remarkable that the algorithm does converge even with $N = 1$.

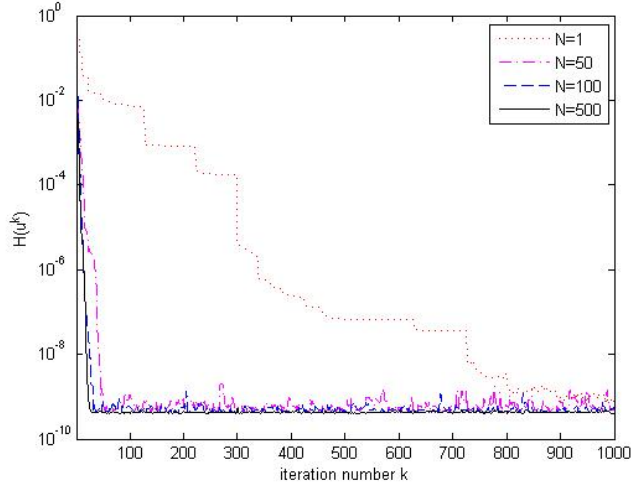

Figure 2. Convergence of mean square error when $N=1$: zoom out view of Figure 1.

Note that when $N=1$, we need to restart the algorithm after every iteration because the covariance matrix has rank 1 and there is no need to search in any directions other than the direction of the residual. In this case, Algorithm SCG effectively becomes the steepest decent method.

## 6.2 Digital Predistortion

Next we present simulation results in applying the SCG_MUL algorithm to digital predistortion of HPA as studied in [1,2]. In [1,2], memory polynomials are used as a predistorter and the computation of the predistorter amounts to solving the minimization Problem 2, in which the basis functions are polynomials. A block diagram of the polynomial predistortion is shown in Figure 3.

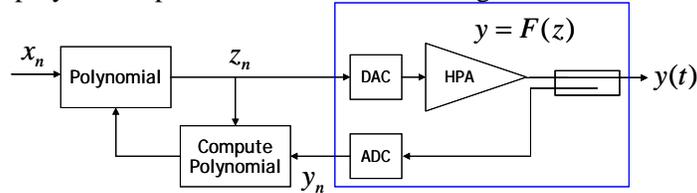

Figure 3. Polynomial predistorter

In a polynomial predistorter, the signal is predistorted by a polynomial. The signal after the predistorter is converted to an analog signal and transmitted to a HPA. A feedback signal from the HPA is sampled, and the pair $y_n, z_n$ forms the input, output pair of the function to be approximated. The objective is to approximate the inverse of the HPA, where $y_n$ is considered the input, and $z_n$ the output. The polynomial that best approximates the inverse function is computed, and it is used as predistorter to form $z_n$ from $x_n$.

In our simulations, we use the memory polynomial PA model as given in Example 2 of [1]. The memory polynomial of degree 5 with 3 delay taps is used as the predistorter. That is, the parameters of Problem 2 are given by $M=5, Q=3$. As suggested by [1], this choice of parameters results in good performance for the predistortion with the given PA model. A multivariate function $g(y)$ is a function of the variables $y_n, y_{n-1}, y_{n-2}$, i.e., $z_n = g(y_n, y_{n-1}, y_{n-2})$.

An OFDM signal with the 16QAM modulation is used in our simulations. At the beginning of each simulation, 25,600 samples are captured for each of $y_n, z_n$. These samples are used to estimate the probability density function $\rho$ using the histogram method. The orthogonal polynomial basis $\{\psi_0(|x|),...,\psi_{M-1}(|x|)\}$ is then formed using the three term recursion. These functions are represented by LUTs with $B = 4096$ entries each. In each data capture of the SCG_MUL algorithm, a total of $N = 1280$ samples are taken for each of $y_n, z_n$.

Three simulations using the SCG_MUL algorithm are performed, which are named Sim1, Sim2 and Sim3. In Sim1, the algorithm is performed with the maximum number of iterations per set of samples captured, i.e., $MQ = 15$ iterations are performed for each data captured. The solution at the last iteration corresponds to the solution by a direct method applied to the normal equations with the inner products formed with the given set of the samples. In Sim2, only one iteration is performed for each set of samples captured. In both Sim1 and Sim2, the weight function used for the orthogonalization of the basis functions is the estimated probability density function. Sim3 is similar to Sim2, but the weight function is the uniform distribution (that is, no weight function is used in the orthogonalization process). In all simulations, the search direction $v$ is reset in every $MQ = 15$ iterations. In Sim1, the reset is performed every time a new set of samples is captured. In all simulations, at each iteration of the SCG_MUL, the newly updated $u$ is immediately used in the outer loop as the predistorter.

In each simulation, a total of 210 SCG iterations are performed. In Sim1, a total of 14 data captures are performed (there are 15 SCG iterations for each data capture), and 210 data captures are performed in other simulations. At the end of each simulation, the linearization of the HPA is achieved, because the HPA output signal $y_n$ is almost identical to the original signal $x_n$. A plot of the spectra of different signals for Sim2 is shown in Figure 4. Spectra from Sim1 and Sim3 are indistinguishable from Figure 4.

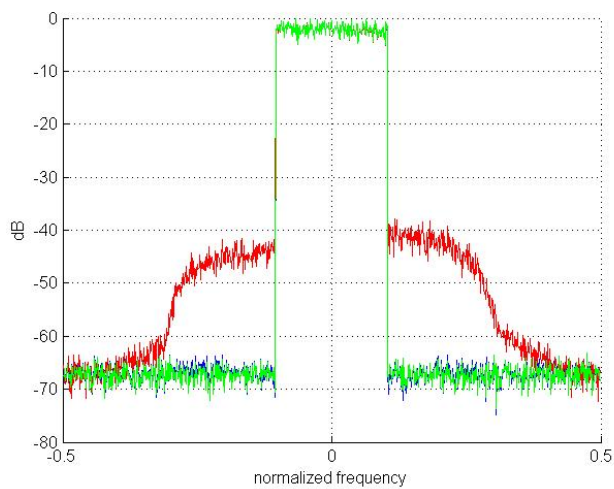

Figure 4. Power spectra

In Figure 4, the red curve, the curve with the highest spectrum shoulders, shows the spectrum of the signal after HPA when no predistortion is used. The green and blue curves, which are almost identical and indistinguishable, show the spectra of the original signal $x_n$ and the signal $y_n$ after HPA when the predistortion is used, respectively.

The spectrum plot in Figure 4 shows that at the end of the iteration after the convergence, the solutions from all three simulations have the desired accuracy because the nonlinear efforts of the HPA have been completely removed with the computed predistorter. A comparison of Figure 4 with the plot in [1, Figure 4] also demonstrates that the solutions from iterative methods of Sim1, Sim2 and Sim3 all have the same accuracy as a solution obtained by a direct method from [1]. However, the SGC algorithm is much more efficient in terms of computational complexity than a direct method.

To show the performance of the SCG_MUL algorithm, we examine the residual computed at the beginning of each SCG iteration. Let $y_n, z_n$ be the set of the captured samples. Let $P$ be the computed polynomial from the previous iteration. Then the normalized residual is defined as $r = \|z_n - P(y_n)\| / \|z_n\|$. The normalized residuals as functions of SCG iteration number are shown in Figure 5 for all simulations.

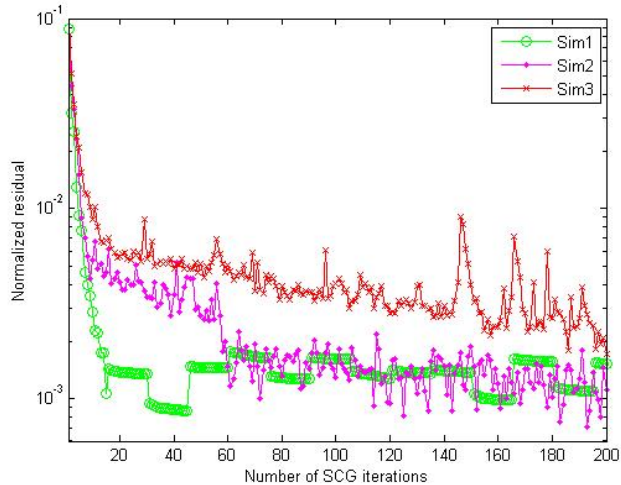

Figure 5. The convergence of normalized residuals

We can make the following observations. First, as expected, in Sim1, a fairly accurate solution is obtained in 15 SCG iterations. After that, the residual does not change significantly. The variation in the residuals after 15 SCG iterations is mainly due to the fact that they are computed with different sets of samples. After convergence, the residual remains almost constant during the 15 SCG iterations in which the same set of samples is used. This demonstrates that after some initial time, there is no need to perform more iterations in the SCG_MUL algorithm using the same set of captured samples. One reason for the residual to be prohibited from being further reduced is because the inverse function of the HPA is not a polynomial. There is also noise in the signal which contributes to errors in the approximation.

In Sim2, the first 5 or 6 iterations are almost identical to those in Sim1. After that, the convergence slows down. Again this is expected because the SCG_MUL algorithm loses orthogonality when different samples are taken at different iterations. However, the normalized residuals are reduced to the similar level as in Sim1 after about 60 SCG iterations (which is equivalent to 4 data captures in Sim1).

In Sim3, we see that the convergence of the SCG_MUL is significantly slower when the uniform distribution is used as the weight function in forming the orthogonal basis. This is because the condition number of the covariance matrix is larger when the weight function is not equal to the probability density function of the samples.

## 7. Conclusions

We have presented a stochastic conjugate gradient method in which the approximating function is computed directly without using the covariance matrix. This reduces complexity of the computation. Furthermore, the inner products involved in the algorithm are computed by evaluating at the samples. We have provided a rigorous proof that the algorithm is convergent in probability. Simulations are performed to confirm the theoretical results.

**Appendix**

**Proof of Lemma 1.**
In iteration $k$ of the SCG method, the update $u^k$ is computed to minimize the functional $J_k(u)$. Therefore, the update satisfies
$$J_k(u^k) \leq J_k(u^{k-1}). \tag{A.1}$$
Combining equations (A.1) and (3.4), we have
$$H_k(u^k) \leq H_k(u^{k-1}). \tag{A.2}$$
If we define
$$\delta_k = \begin{cases} H_k(u^k)/H_k(u^{k-1}), & \text{if } H_k(u^{k-1}) \neq 0 \\ 0, & \text{otherwise} \end{cases}, \tag{A.3}$$
then (3.10) holds, and equations (A.2) and (A.3) imply (3.9).

In order to show (3.11), we note that at the start or the reset of Algorithm SCG, the computation is exactly the same as in the steepest decent method because the search direction is the residual. By using the same process as in the steepest decent method, see for example [1], it can be easily shown that
$$H_j(u^j) \leq (1 - c_j^{-1})H_j(u^{j-1}), \tag{A.4}$$
if iteration $j$ is at a reset.

Algorithm SCG requires a reset at least once in every $M$ iterations, and therefore, (A.4) is satisfied at least once in every $M$ iterations. For large $k$, there must be at least $k/M$ occurrences of reset, and hence, as many times for which (A.4) is satisfied. For each such occurrence, we have
$$\delta_j \leq 1 - c_j^{-1} \leq 1 - c_0^{-1}. \tag{A.5}$$
In arriving at (A.5), assumption (3.8) has been applied. This shows that there is a fraction of $k$ iterations for which (A.5) is satisfied. More precisely, there exists a $p$ with $\frac{1}{M} \leq p \leq 1$, so that for large $k$, among all $\delta_j$, $j = 1,...,k$, at least $pk$ of them satisfy (A.5), and the others satisfy (3.9). This proves (3.11).

**Proof of Lemma 2.**
Since $H(u^k) \neq 0$ for all $k = 1,2,...$, we define
$$\begin{aligned}\omega_{k0} &= H_k(u^k)/H(u^k), \\ \omega_{k1} &= H_k(u^{k-1})/H(u^{k-1}).\end{aligned} \tag{A.6}$$
From (3.10) and (A.6), we have
$$\omega_{k0} H(u^k) = \delta_k \omega_{k1} H(u^{k-1}). \tag{A.7}$$
We first give the proof by assuming $\omega_{k0} \neq 0$ and $\omega_{k1} \neq 0$ for all $k$. Applying (A.7) repeatedly for $k, k-1,...,1$, we get
$$\begin{aligned}H(u^k) &= \delta_k \omega_{k1} \omega_{k0}^{-1} H(u^{k-1}) = \prod_{j=1}^{k}\delta_j \prod_{j=1}^{k}\omega_{j1}\left(\prod_{j=1}^{k}\omega_{j0}\right)^{-1} H(u^0) = \eta_k \theta_k H(u^0), \\ \eta_k &= \prod_{j=1}^{k}\delta_j, \\ \theta_k &= \prod_{j=1}^{k}\omega_{j1}\left(\prod_{j=1}^{k}\omega_{j0}\right)^{-1}.\end{aligned} \tag{A.8}$$
Therefore, $\eta_k$ satisfies (3.13) with $p_0 = p$ thanks to Lemma 1. The quantities $\omega_{j0}, \omega_{j1}$ are random variables defined in (3.6). Since they are independent and identically distributed by assumption, the product $\prod_{j=1}^{k}\omega_{j1}$ approaches to a log-normal distribution as $k \to +\infty$ according to Central Limit

Theorem. Since the mean and variance of the logarithm of the random variables in (A.6) are $\mu$ and $\sigma^2/2$, respectively, the log-normal distribution in the limit has the parameters $k\mu$ and $k\sigma^2/2$, respectively. The same argument also applies to product $\prod_{j=1}^{k} \omega_{j0}$. Furthermore, the inverse of a log-normal distribution is also a log-normal distribution, and the product of independent log-normal distributions is also log-normal. Therefore, for large $k$, we have the following distributions

$$\prod_{j=1}^{k} \omega_{jq} \sim \text{LogN}(k\mu, k\sigma^2/2), \ q = 0, 1$$

$$\left(\prod_{j=1}^{k} \omega_{j0}\right)^{-1} \sim \text{LogN}(-k\mu, k\sigma^2/2),$$

$$\theta_k = \prod_{j=1}^{k} \omega_{j1} \left(\prod_{j=1}^{k} \omega_{j0}\right)^{-1} \sim \text{LogN}(0, k\sigma^2).$$

The probability density function of the log-normal distribution $\text{LogN}(0, k\sigma^2)$ is given by (3.14). This proves Lemma 2 under the assumption of $\omega_{k0} \neq 0$ and $\omega_{k1} \neq 0$ for all $k$.

Next, we consider the case where $\omega_{j0}$ or $\omega_{j1} = 0$ for some $j = 1,...,k$. Because of (A.7) and the assumption that $H(u^j) \neq 0$, we have $\omega_{j1} = 0$ implies $\omega_{j0} = 0$. Therefore, we only need to consider the case where $\omega_{j0} = 0$ for some $j$. From (A.6), $\omega_{j0} = 0$ implies $H_j(u^j) = 0$, and hence $u^j = \bar{u}^j$. For each of such $j$, we have, because of (3.8),

$$H(u^j) = \langle \bar{u}^j - \bar{u}, \bar{u}^j - \bar{u} \rangle \leq \left(\sqrt{\langle \bar{u}^j, \bar{u}^j \rangle} + \sqrt{\langle \bar{u}, \bar{u} \rangle}\right)^2 \leq \left(d_0 + \sqrt{\langle \bar{u}, \bar{u} \rangle}\right)^2 \stackrel{\Delta}{=} d_1 H(u^0). \quad (A.9)$$

This shows that (A.8) is still true if we replace $\delta_j \omega_{j1} \omega_{j0}^{-1}$ by $d_1$ in each of the terms in (A.8) for which $\omega_{j0} = 0$. Since $u^j$ is computed by random samples and only one iteration of the SCG is performed, the probability of $\omega_{j0} = 0$, hence $H_j(u^j) = 0$, is zero, i.e., $\Pr(H_j(u^j) = 0) = 0$. Therefore, for large $k$, there are no more than, say, $(p/4)k$ terms in (A.8) for which (A.9) holds, and for the rest of terms (A.7) holds. Absorbing those terms for which (A.9) holds into the product $\eta_k$ of (A.8), we have

$$\eta_k \leq \left(1 - c_0^{-1}\right)^{3pk/4} d_1^{pk/4} \leq \left(1 - c_0^{-1}\right)^{pk/4}, \text{ for large } k.$$

Therefore, (3.13) holds for some $p_0$, say $p_0 = p/4$. The product $\theta_k$ has few terms than given in (A.8), but for large $k$, it still approximates the log-normal distribution with the same parameters, which concludes proof.

**Proof of Theorem 1.**
First, if $H(u^k) = 0$ for some $k$, then $\langle u^k - \bar{u}, u^k - \bar{u} \rangle = 0$. This shows $u^k = \bar{u}$, and $u^k$ is therefore the solution to Problem 1. Otherwise, if $H(u^k) \neq 0$ for all $k = 1, 2,...$, then Lemma 2 holds. For a given $\varepsilon > 0$, we consider the probability of $|H(u^k)| < \varepsilon$. From Lemma 1 and Lemma 2, we have

$$|H(u^k)| = \eta_k \theta_k H(u^0) \leq \left(1 - c_0^{-1}\right)^{p_0 k} H(u^0) \theta_k. \quad (A.10)$$

Since LHS of (A.10) is smaller than or equal RHS, the probability of LHS being smaller than a number is larger than the probability of RHS being smaller than the same number. Therefore, we have

$$\Pr\left(|H(u^k)| \leq \varepsilon\right) \geq \Pr\left(\left(1 - c_0^{-1}\right)^{p_0 k} H(u^0) \theta_k \leq \varepsilon\right). \quad (A.11)$$

The probability on the RHS of (A.11) is the same as that of $\theta_k \leq \varepsilon_k$ if we define

$$\varepsilon_k = \varepsilon\left(1 - c_0^{-1}\right)^{-p_0 k}\left(H(u^0)\right)^{-1}, \tag{A.12}$$

because all numbers on the RHS of (A.12) are deterministic positive numbers. Therefore,

$$\Pr\left(|H(u^k)| \leq \varepsilon\right) \geq \Pr\left(\theta_k \leq \varepsilon_k\right). \tag{A.13}$$

Since the probability density function $\theta_k$ is given by (3.14), we have

$$\Pr\left(\theta_k \leq \varepsilon_k\right) = \int_0^{\varepsilon_k} f_{\theta_k}(t)dt = \int_0^{\varepsilon_k} f_{\text{LogN}}(t;0,k\sigma^2)dt. \tag{A.14}$$

Thus, combing (A.13), (A.14), we get

$$\Pr\left(|H(u^k)| \leq \varepsilon\right) \geq \int_0^{\varepsilon_k} f_{\text{LogN}}(t;0,k\sigma^2)dt. \tag{A.15}$$

To evaluate the integral in (A.15), we make a change of variable $s^{\sigma\sqrt{k}} = t$. Then the integral is computed as

$$\int_0^{\varepsilon_k} f_{\text{LogN}}(t;0,k\sigma^2)dt = \int_0^{\varepsilon_k} \frac{1}{t\sqrt{k}\sigma\sqrt{2\pi}} e^{-\frac{(\ln t)^2}{2k\sigma^2}} dt = \int_0^{s_k} \frac{1}{s\sqrt{2\pi}} e^{-\frac{(\ln s)^2}{2}} ds = \int_0^{s_k} f_{\text{LogN}}(s;0,1)ds, \tag{A.16}$$

where

$$s_k = (\varepsilon_k)^{\frac{1}{\sigma\sqrt{k}}} = \varepsilon^{\frac{1}{\sigma\sqrt{k}}} \left(H(u^0)\right)^{-\frac{1}{\sigma\sqrt{k}}} \left(1 - c_0^{-1}\right)^{-p\frac{\sqrt{k}}{\sigma}}. \tag{A.17}$$

As $k \to +\infty$, the first two terms in the RHS of (A.17) goes to 1 and the last term goes to $+\infty$. Therefore, the limit of (A.17) is $\lim_{k \to +\infty} s_k = +\infty$. Now taking the limit of (A.16), we have

$$\lim_{k \to +\infty} \int_0^{\varepsilon_k} f_{\text{LogN}}(t;0,k\sigma^2)dt = \lim_{k \to +\infty} \int_0^{s_k} f_{\text{LogN}}(s;0,1)ds = \int_0^{+\infty} f_{\text{LogN}}(s;0,1)ds = 1. \tag{A.18}$$

Finally, from (A.15) and (A.18), we have

$$\lim_{k \to +\infty} \Pr\left(|H(u^k)| \leq \varepsilon\right) \geq \lim_{k \to +\infty} \int_0^{\varepsilon_k} f_{\text{LogN}}(t;0,k\sigma^2)dt = 1.$$

Since the probability cannot exceed 1, this implies (3.15), which concludes the proof.